\documentclass{article}
\usepackage[latin1]{inputenc}
\usepackage{amsmath}
\usepackage{amsfonts}
\usepackage{amssymb}
\usepackage{latexsym}
\usepackage{cite}
\usepackage[hidelinks]{hyperref}

 \DeclareMathOperator{\tr}{tr}

\newtheorem{te}{Theorem}[section]

\newtheorem{lemma}{Lemma}[section]

\newcommand{\rr}{\mathbb{R}}

\newcommand{\al}{\alpha}
\newcommand{\alp}{\alpha^{*}}

\newcommand{\la}{\lambda}
\newcommand{\rn}{\mathbb{R}^{N}}
\newcommand{\sm}{\mathcal{S}_{N}}

\begin{document}

\title{On a family of  exact solutions for a nonlinear diffusion equation}
\author{Rodrigo Meneses Pacheco\\
\\
\small Escuela de Ingenier\'ia Civil, Facultad de Ingenier\'ia, Universidad de Valpara\'iso\\
\small Avda. Errazuriz 1834, Valpara\'iso, Chile\\
\href{mailto:rodrigo.meneses@uv.cl}{rodrigo.meneses@uv.cl}}

\date{August 11, 2014}
\maketitle
\begin{abstract}
We present a complete description of the similarity solutions $u_{\al}(x,t)=t^{-\alpha/2}f(\Vert x  \Vert/\sqrt{t};\al)$  for  the following nonlinear diffusion  equation
$$ u_{t}+\gamma\vert u_{t} \vert =\Delta u\qquad(-1<\gamma<1) $$
The behaviors of these  solutions are obtained through the explicit representation of  $f(\eta;\al)$,  in terms of  Kummer and Tricomi functions. Considering results about confluent hypergeometric functions, new methods to describe asymptotic and oscillatory
behaviors of the similarity solutions are obtained. We prove that there exists an increasing and unbounded sequence of positive similarity exponents such that the associated profile $f$ has a gaussian rate decay. These special similarity exponents are related with the zeros of  Kummer and Tricomi functions.  Finally, we indicate how to extend  our results  on more
general nonlinear diffusion equations.

bf{keywords:} Similarity solutions, Kummer functions, Tricomi functions, nonlinear diffusion equations
\end{abstract}

\section{Introduction}
 In physics and engineering problems finding similarity solutions is common. This feature is a consequence of symmetries in the underlying models used to solve these problems.  It is known  that these special functions  are frequently used to describe  the  intermediate asymptotics of general solutions, i.e.,  behaviors  in regions where the dependence with initial and/or boundary conditions already disappears \cite{barenblatt1996scaling,barenblatt1972self}.\\
 For nonlinear PDEs, the similarity solutions are also used to study the development of singularities and other class of particular phenomena such as  blow-up or extinction in finite time  \cite{eggers2009role}. Thus, there is an extensive literature on the relation between similarity solutions and behaviors of general solutions  and, therefore, the analysis and description of these special class of solutions have a practical relevant importance from different points of view.

 In this paper, we study the following nonlinear diffusion equation
 \begin{equation}
\label{barenblattcanonico}
u_{t}+\gamma\vert  u_{t} \vert=\Delta u\qquad\textrm{in}\ \rn\times]0,\infty[,
\end{equation}
where $-1<\gamma<1$. This equation was formulated in \cite{barenblatt1955concerning} as a  nonlinear filtration theory  model and since it's presentation it has been studied different framework in applied mathematics \cite{barenblatt1996scaling,khuzhaerov2007inverse,peng2007g} and considered  as a simple model of nonlinear parabolic equation \cite{aronson1994calculation,caffarelli2008counterexample,HUANG2012LARGE,kamin1991barenblatt}.

Our main objective is to obtain, through  explicit representations,  asymptotic and oscillatory behaviors of the similarity solutions
\begin{equation}
\label{selfsimilar}
u_{\al}(x,t)=t^{-\alpha/2}f(\eta;\al),\quad \textrm{with}\quad  \eta= \Vert x  \Vert/t^{1/2}\ \textrm{and}\ \al>0.
\end{equation}
  The analysis is centered on the profiles $f(\eta;\al)$ which are studied  as solutions of a specific nonlinear Cauchy problem obtained when (\ref{selfsimilar}) is substituted into (\ref{barenblattcanonico}). The Cauchy problem to $f(\eta;\al)$ is described by a nonlinear ODE defined by a piecewise linear discontinuous function. Thus, the behaviors of $f(\eta;\al)$ are studied considering a collection of linear Cauchy  problems. \\ Following  ideas from \cite{barenblatt1969self,caffarelli2008counterexample,Kerchman1971self,kroger1994regularity}, we obtain  a closed representation for (\ref{selfsimilar}) in terms of confluent hypergeometric functions, and oscillatory and asymptotic behaviors are obtained using the asymptotic representations of Kummer and Tricomi functions.

Our first main result reads as follows.
\begin{te}
\label{principal1}
There exists an increasing and unbounded sequence of similarity exponents
\begin{equation}
\label{valorespropiosint}
0<\alpha_{0}<\alpha_{1}<\alpha_{2}<\dots
\end{equation}
such that $f(\eta;\al_{k})\sim \eta^{\al_{k}-N}e^{-\frac{(1+(-1)^{k}\gamma)}{4}\eta^{2}}$  for $\eta$ large. When  $\al\neq\al_{k}$, the asymptotic behavior is given by $f(\eta;\alpha)\sim \eta^{-\alpha}$.
\end{te}
We show that the exponents  in (\ref{valorespropiosint}) are related with zeros of confluent hypergeometric functions.
Specifically, the exponent $\al_{m}$ in (\ref{valorespropiosint}) are determined by the analysis of the following system
\begin{equation}
\label{sistemaexponente}
\begin{array}{rcl}
M(\frac{N-\al}{2}-1,\frac{N}{2};\frac{(1-\gamma)}{4}\eta_{1}^{2})&=&0\\
\\
U(\frac{N-\al}{2}-1,\frac{N}{2};\frac{(1+(-1)^{m}\gamma)}{4}\eta^{2}_{m+1})&=&0\\
\\
M(\frac{N-\al}{2}-1,\frac{N}{2};\frac{(1+(-1)^{m}\gamma)}{4}\eta_{m+1}^{2}) \cdot \Gamma(\lbrace N-\al \rbrace/2)\cdot (-1)^{m}&>&0
\end{array}
\end{equation}
obtaining a complete generalization  of the results presented in \cite{barenblatt1969self,Kerchman1971self} for the exponent $\al_{0}$  for $N=1$ and $N$ arbitrary, respectively.\\
Assuming that there exists exactly $m+1$ roots $\eta_{1}<\eta_{2}<\dots<\eta_{m}<\eta_{m+1}$ for $\eta f'+\al f=0$, in (\ref{sistemaexponente}), the values $\eta_{1}$ and $\eta_{m+1}$ denotes the smallest and largest roots, respectively.  We use $M(a,b;z)$, $U(a,b;z)$ to denote the Kummer and Tricomi functions  and $\Gamma$ is the symbol for the Gamma-Euler function.

 To complete the analysis about the behaviors of $f(\eta;\al)$,  in our second main result the following oscillatory behaviors are presented.
\begin{te}
\label{principal2}
Let $\al_{k}$ be a similarity exponent in Theorem \ref{principal1}. For each $\al$ such that $\al_{k}<\al\leq\al_{k+1}$,
the profiles $f(\eta;\al)$ have exactly $k$ zeros in $[0,\infty[$. If $\al\leq \al_{0}$ follows $f(\eta;\al)>0$.
\end{te}

This paper is organized as follows: In Section \ref{preliminares} we present results about (\ref{barenblattcanonico}) and a sequence of lemmas necessary for our proof on the oscillatory behavior of $f(\eta;\al)$.
The method for the representation of similarity solutions in terms of Kummer and Tricomi functions is presented in  Section \ref{sectionsimilarity}. Theorems \ref{principal1} and \ref{principal2} are proved in Section \ref{principal}, and the relation between
the spatial decays of $f$ and the zeros of Kummer and Tricomi functions is showed.
Finally, in Section \ref{extension}, we present how  Theorems \ref{principal1} and
\ref{principal2} can be extended to more general nonlinear diffusion equations.

\section{Preliminaries}\label{preliminares}

\subsection{About the Barenblatt equation}

Equation (\ref{barenblattcanonico}), commonly called Barenblatt equation, was formulated  in  \cite{barenblatt1955concerning} as a simple model to study the filtration of a low compressibility fluid through porous elasto-plastic media where deformations are irreversible. Such behavior is modeled by considering a piezo-conductivity coefficient which can take different values.Therefore, according to Darcy's law the following is observed
\begin{equation}
\label{barenblatt}
p_{t}=\kappa(p_{t})\Delta p,\qquad \kappa(z)=\left\lbrace \begin{array}{ll}\kappa^{+}& (z>0)\\ \kappa^{-}& (z\leq0)  \end{array}  \right.
\end{equation}
where $\kappa^{-}$ and $\kappa^{+}$ depends upon the physical parameters of the model such as: coefficient of permeability; compressibility of the porous medium; compressibility of the fluid and others. \\
Considering the following transformations:
\begin{equation*}
\label{defpar}
 \gamma=\frac{\kappa^{-}-\kappa^{+}}{\kappa^{-}+\kappa^{+}}\quad\textrm{and}\quad t=\frac{1}{2}\left(\frac{1}{\kappa^{+}} -\frac{1}{\kappa^{-}}  \right)t^{*},
\end{equation*}
equation (\ref{barenblatt}) can be written as (\ref{barenblattcanonico}), where $u(x,t^{*})=p(x,t)$ (omitting the asterisk in (\ref{barenblattcanonico})), satisfying directly the condition $-1<\gamma<1$.  For further details about the derivation  and physical consideration of this equation see \cite{barenblatt1996scaling,barenblatt1989theory,barenblatt1955concerning} and their references.

Recently, equation (\ref{barenblatt})  was considered as a model  in sublinear expectation theory  where the notions of G-Normal distribution, G-expectation and G-Brownian motion are introduced through the nonlinear parabolic equation
\begin{equation}
\label{barenblatt2}
u_{t}-\frac{1}{2}((u_{xx})^{+}-\sigma^{2} (u_{xx})^{-})=0;\qquad 0\leq \sigma\leq 1
\end{equation}
with  $a^{+}=\max\lbrace 0,a\rbrace$, $a^{-}=(-a)^{+}$, i.e. $\gamma<0$ in (\ref{barenblattcanonico}). In the framework of  sublinear expectation theory, equation (\ref{barenblatt2}) is often called   G-Heat equation. More details about its application can be studied in \cite{hu2012explicit, peng2007g, peng2008multi}.

Due their simple mathematical structure, equation (\ref{barenblattcanonico}) often has been considered   as practical example of nonlinear diffusion equations where is possible to use different tools for the analysis. The mathematical analyzes on (\ref{barenblattcanonico}) are related principally with the similarity solution associated with $\al=\al_{0}$, so-called anomalous exponent \cite{aronson1994calculation,goldenfeld1990anomalous} .\\
Recently, Eq. (\ref{barenblattcanonico}) was used to study  qualitative results related with more general diffusion equations, specifically, for the study of large time geometrical properties \cite{HUANG2012LARGE} and for the study of a counterexample to regularity for fully nonlinear equations \cite{caffarelli2008counterexample}.\\
Finally,  we remark that in \cite{hu2012explicit} an explicit family of similarity solution of (\ref{barenblatt2})  with polynomial initial condition $u(x,0)=x^{n}$ was presented. The author  use  a matching between linear problems, similar to our method to solve (\ref{barenblattcanonico}).

\subsection{Nonlinear ODE for the profiles $f(\eta;\al)$}

If we assume that the function in (\ref{selfsimilar}) is a classical solution of (\ref{barenblattcanonico}), when we consider $\eta>0$  such that $( \alpha f+ \eta f  )>0$, the following ODE is obtained:
\begin{equation*}
\label{rama1}
f''+\frac{N-1}{\eta}f'=-\frac{(1-\gamma)}{2} ( \alpha f+ \eta f  ),
\end{equation*}
while in cases $( \alpha f+ \eta f  )<0$ follows
\begin{equation*}
\label{rama2}
f''+\frac{N-1}{\eta}f'=\frac{(1+\gamma)}{2} ( \alpha f+ \eta f  ).
\end{equation*}
We note that in the formulation of the ODE  we need the following functions:
\begin{equation}
\label{fundiscontinua}
\sigma(z)=\left\lbrace \begin{array}{ll}-(1-\gamma)/2& (z>0)\\ -(1+\gamma)/2& (z\leq0)  \end{array}  \right.
\end{equation}
On the other hand, for the regularity of  $u_{\al}(x,t)$  at  $x=0$, it's necessary  that $f'(0)=0$.
Notice that, when multiplying (\ref{barenblattcanonico}) by a positive constant $A>0$, we get the solution $Au_{\al}(x,t)$ and therefore, the value $f(0)=A>0$ can be chosen arbitrarily. The profile $f(\eta)$ is conveniently normalized by the relation $f(0)=1$. \\
Thus, the profile $f(\eta;\al)$   satisfies the following nonlinear Cauchy problem
\begin{equation}
\label{cauchynolinear}
\left\lbrace  \begin{array}{rcl}
                f''+\frac{N-1}{\eta}f'&=&\sigma(\alpha f+ \eta f' ) (  \alpha f+ \eta f'  )\\
                                 f(0) &=&1\\
                                 f'(0)&=&0
              \end{array}
\right.
\end{equation}
 The existence and uniqueness  follows from the fact that $\sigma(z)$ is a piecewise linear function. Thus, for each $\al>0$ there exist an unique solution for (\ref{cauchynolinear}) and using these profiles we obtain a one-parameter family of similarity solutions (\ref{selfsimilar}) for the equation (\ref{barenblattcanonico}).

From the form of the ODE in (\ref{cauchynolinear}), the solution  can be obtained using a collection of linear Cauchy problems, each one formulated in different ranges of $\eta$. Each of these ranges correspond to an interval where the term $\eta f'(\eta;\al)+\al f(\eta;\al)$
does not change sign. Hence, considering the $n$ positive roots $\eta_{1}<\eta_{2}<\dots<\eta_{n}$ of
\begin{equation}
\label{ecuacionparalosceros}
\eta f'(\eta;\al)+\al f(\eta;\al)=0,
\end{equation}
the ranges of linear systems are: $]0,\eta_{0}[,\ ]\eta_{0},\eta_{1}[,\dots,\ ]\eta_{n-1},\eta_{n}[,\  ]\eta_{n},\infty[$ (see (\ref{formak}) below).\\
Through an appropriate transformation, the ODE in each linear problem is associated with
the following Kummer ODE:
  \begin{equation}
  \label{kummereq}
  zu''+((N/2)-z)u'-(\al/2)u=0.
\end{equation}
Hence, the profile $f$ can be represented through of the solutions of (\ref{kummereq}). Eq. (\ref{kummereq})   has an entire solution called Kummer function, defined by
\begin{equation}
\label{kummer}
M(\al/2,N/2;z)=1+\sum_{k=1}^{\infty}\frac{(\al/2)_{k}}{(N/2)_{k}}\frac{z^{k}}{k!},
\end{equation}
where $(\la)_{k}$ is the Pochhammer symbol:
\begin{equation}
\label{simbolo}
(\la)_{0}=1,\quad (\la)_{k}=\frac{\Gamma(\la +k)}{\Gamma(k)}=\la(\la+1)(\la+2)\dots(\la+k-1),\qquad k=1,2,\dots ,
\end{equation}
The Kummer function is also called confluent hypergeometric function (of the first kind) and it's also denoted by $_{1}F_{1}(a,b;z)$.

We remark that  when $\al< N-2$ follows $\al f+\eta f'>0$ and therefore $f$ is obtained solving (\ref{rama1}). Thus, if $\al<N-2$, we have the following representation
\begin{equation}
\label{repre1}
f(\eta;\al)=M(\frac{\al}{2},\frac{N}{2},-\frac{(1-\gamma)}{4}\eta^{2}),
\end{equation}
satisfying $f(\eta;\al)>0$. This result is obtained  considering $z=-(1-\gamma)\eta^{2}/4$ and the transformation  $g(z)=f(\eta;\al)$. We remark that the existence of roots for $\eta f'+\al f=0$ for each $\al>N-2$  is studied in Lemma \ref{cambiodesigno1} (below).

On the other hand, in the linear case  $\gamma=0$,   the family of self-similar solutions with $\al>0$  are described through the Kummer functions (\ref{kummer}). Specifically, if $\gamma=0$   the profile $f(\eta;\al)$ is described as follows:
\begin{equation*}
\label{solheat}
f_{heat}(\eta;\al)=M(\frac{\al}{2},\frac{N}{2};-\frac{\eta^{2}}{4}).
\end{equation*}
In cases  $b-a=-l$ $(l=0,1,2,\dots)$, the Kummer function $M(a,b;z)$ has an
exponential-type behavior to $z\to-\infty$. This result is directly verified using the Kummer transformation
\begin{equation}
\label{decaimiento}
\begin{array}{rcl}
M(a,b;z)&=&e^{-z}M(b-a,b;-z)\\
        &=&e^{-z}\left\lbrace\displaystyle{ 1+\sum_{k=0}^{\infty} \frac{(b-a)_{k}}{(b)_{k}}\frac{(-z)^{k}}{k!}} \right\rbrace.
\end{array}
\end{equation}
Now, from the definition of $(a)_{k}$ in (\ref{simbolo}), when $b-a=-l$ (with $l=0,1,2,\dots$), we have $(b-a)_{k}=0$ for each $k\geq l$ and therefore we obtain de exponential behavior. Using the representation (\ref{decaimiento}), we note that $f_{heat}(\eta;\al)$ have an
exponential-type decay when $\frac{\al-N}{2}=l$ $(l=0,1,2,\dots)$. Hence, when $\gamma=0$ the similarity exponents $\al_{l}$ in (\ref{valorespropiosint}) are represented by $\al_{l}=N+2l$.

In the case where $\gamma\neq 0$, we need to use recurrence results and asymptotic representations of the  Tricomi function or confluent hypergeometric function
of the second kind, denoted by $U(\al/2,N/2;z)$, see for instance \cite{abramowitz1972handbook,olver1974introduction}.

\subsection{Outline of the construction of the Similarity Solutions}
We continue using $f(\eta;\al)$ to denote the solution of  (\ref{cauchynolinear}). If (\ref{ecuacionparalosceros}) has no solutions the representation for $f$ is given in (\ref{repre1}). In other case, we need solve more linear Cauchy problems. \\
Let $\eta_{1}>0$ be the first positive root of (\ref{ecuacionparalosceros}). Using the initial conditions for $f$, we have $\al f+\eta f'>0$ when $\eta=0$. Therefore for each $\eta$  such that $0<\eta<\eta_{1}$ follows $\al f+\eta f'>0$ . Let $f_{1}(\eta;\al)$ be  the solution of the following linear problem:
\begin{equation}
\label{cauchylinear1}
\left\lbrace  \begin{array}{rcl}
                f_{1}''+\frac{N-1}{\eta}f_{1}'&=&-\frac{(1-\gamma)}{2}(  \alpha f_{1}+ \eta f_{1}'  )\\
                                 f_{1}(0) &=&1\\
                                 f_{1}'(0)&=&0
              \end{array}
\right.
\end{equation}
As $\al f+\eta f'>0$ when $0<\eta<\eta_{1}$, from $\sigma(z)$ in
(\ref{fundiscontinua}), we have
$$f''+\frac{N-1}{\eta}f'=-\frac{(1-\gamma)}{2}(  \alpha f+ \eta f'  )\qquad (0<\eta<\eta_{1})$$
and thus, through classical uniqueness results on ODE, follows
\begin{equation}
\label{1comp}
 f(\eta;\alpha)=f_{1}(\eta;\alpha),\qquad \textrm{para}\ 0\leq\eta\leq\eta_{1}.
\end{equation}
Now, we continue in a similar way. Assume that (\ref{ecuacionparalosceros}) has exactly $n$ positive roots, denoted by $\eta_{1}<\eta_{2}<\dots<\eta_{n}$. Through these roots, we describe different ranges where auxiliary  linear Cauchy  problems are
defined.\\
Let $f_{2}(\eta;\al)$ be  the solution of the following Cauchy problem
   \begin{equation}
\label{cauchylinear2}
\left\lbrace  \begin{array}{rcl}
                f_{2}''+\frac{N-1}{\eta}f_{2}'&=&-\frac{(1+\gamma)}{2}(  \alpha f_{2}+ \eta f_{2}'  )\qquad\eta>\eta_{1}\\
                                 f_{2}(\eta_{1}) &=& f_{1}(\eta_{1})\\
                                 f_{2}'(\eta_{1})&=& f_{1}'(\eta_{1})
              \end{array}
\right.
\end{equation}
From the definition of $\eta_{i}$, follows $\al f+\eta f'<0$ in $\eta_{1}\leq \eta\leq \eta_{2}$. Hence,
$$f''+\frac{N-1}{\eta}f'=-\frac{(1+\gamma)}{2}(  \alpha f+ \eta f'  )\qquad (\eta_{1}<\eta<\eta_{2})$$
with $f_{2}(\eta_{1}) = f(\eta_{1})$ and $f_{2}'(\eta_{1})= f'(\eta_{1})$. Using  uniqueness results, we have
\begin{equation*}
\label{2comp}
 f(\eta;\alpha)=f_{2}(\eta;\alpha),\qquad  \eta_{1}\leq\eta\leq\eta_{2}.
\end{equation*}
Following the previous ideas, the profile $f(\eta;\al)$ is defined as piecewise
function
\begin{equation}
\label{formak}
f(\eta;\al)=\left\lbrace    \begin{array}{ccl}
                             f_{1}(\eta;\al)  & & 0\leq \eta\leq \eta_{1}\\
                             f_{2}(\eta;\al)  & & \eta_{1}\leq  \eta\leq \eta_{2}\\
                             \vdots\\
                             f_{n}(\eta;\al)  & & \eta_{n-1}\leq  \eta\leq \eta_{n}\\
                             f_{n+1}(\eta;\al)& & \eta_{n}\leq  \eta\\
                            \end{array}
             \right.
\end{equation}
where $f_{m}(\eta;\al)$ is the solution of the linear Cauchy problem:
   \begin{equation*}
\label{cauchylineari}
\left\lbrace  \begin{array}{rcl}
                f_{m}''+\frac{N-1}{\eta}f_{m}'&=&-\frac{(1+(-1)^{m}\gamma)}{2}(  \alpha f_{m}+ \eta f_{m}'  )\qquad\eta_{m-1}<\eta<\eta_{m}\\
                                f_{m}(\eta_{m-1}) &=& f_{m-1}(\eta_{m-1})\\
                                 f_{m}'(\eta_{m-1})&=& f_{m-1}'(\eta_{m-1})
              \end{array}
\right.
\end{equation*}
The case $m=1$ is given in (\ref{cauchylinear1}) and in the case $m=n+1$ we consider $\eta_{n+1}=\infty$ in (\ref{cauchylineari}).

Finally, for each $m=1,2,\dots, n+1$, we consider the following change of variable and transformation:
\begin{equation}
 \label{transformacion}
 z_{m}=-\frac{(1+(-1)^{m}\gamma)}{4}\eta^{2};\qquad g_{m}(z_{m};\alpha)=f_{m}(\eta;\alpha)
 \end{equation}
Considering these substitutions,  the ODE in (\ref{cauchylineari}) is rewritten as follows:
 \begin{equation}
 \label{eqk}
 z_{m}g_{m}''+((N/2)-z_{m})g_{m}'-(\alpha/2)g_{m}=0,
 \end{equation}
 i.e, a Kummer-type ODE with parameters $a=\al/2$ and $b=N/2$. Thus, the functions $f_{m}(\eta;\al)$ can be represented by the Kummer and Tricomi  functions. For $m=1$ we use (\ref{repre1}) to represent $f_{1}(\eta;\al)$. Hence, knowing (\ref{1comp}), from (\ref{repre1}) we obtain the representation of $f(\eta;\al)$ when $\eta\in[0,\eta_{1}]$. In Section \ref{sectionsimilarity} we return with the representation of similarity solutions through the confluent hypergeometric functions.

In regard to the matching technique between linear problems, in  \cite{caffarelli2008counterexample,kroger1994regularity}  this technique was used to study similarity solutions of  nonlinear problems related to (\ref{barenblattcanonico}). In  \cite{kroger1994regularity} an equation with similar symmetries  as (\ref{barenblattcanonico})  was  studied (see (\ref{extremalpucci}) below). The author uses the behaviors of confluent hypergeometric functions to obtain the existence of an exponent such that the similarity solution vanishes at least exponentially fast at infinity.\\
On the other hand, in \cite{caffarelli2008counterexample}  a method to describe similarity solutions of a nonlinear diffusion problem was presented to obtain  a counterexample for the regularity in nonlinear diffusion problems. In this article the authors showed a similar  results  as presented in Theorems \ref{principal1} and \ref{principal2} in our paper, but the results on  similarity solutions are not obtained through the explicit representations.

 To finish this section we present  the following lemma where a sufficient condition for the existence of $\eta_{1}<\infty$ is presented.
   \begin{lemma}
 \label{cambiodesigno1}
 The equation $\al f+\eta f'=0$ has solutions if and only if $\al>N-2$.
 \end{lemma}
 \noindent{\bf{Proof:}} Using the representation (\ref{1comp}) we know that the first
 positive solution of (\ref{ecuacionparalosceros}) is related with the equation
 $$ \eta f_{1}'(\eta;\al)+\al f_{1}(\eta;\al)=0  $$
 where  $f_{1}(\eta;\al)$ is  solution of the linear problem (\ref{cauchylinear1}).  Considering $m=1$ in (\ref{transformacion}), we get (\ref{kummereq}) and therefore the solution is represented as follows:
\begin{equation*}
\label{representacion1}
f_{1}(\eta;\al)=M(\frac{\al}{2},\frac{N}{2};-\frac{(1-\gamma)}{4}\eta^{2}),
\end{equation*}
with $M(a,b;z)$  defined in (\ref{kummer}). Thus, using this representation, we notice that the root
$\eta_{1}$ is determined by the following equation
\begin{equation}
\label{primeraraiz}
\al M(\frac{\al}{2},\frac{N}{2};z) +2z  M'(\frac{\al}{2},\frac{N}{2};z) =0,
\end{equation}
with  $z=-(1-\gamma)\eta^{2}/4$. Through the recurrence relation (see 13.4.10 in \cite{abramowitz1972handbook})
\begin{equation}
\label{recurrenciaKummer}
zM'(a,b;z)+aM(a,b;z)=aM(a+1,b;z),
\end{equation}
equation (\ref{primeraraiz}) is rewritten as follows:
$$ M(\frac{\al}{2}+1,\frac{N}{2};-\frac{(1-\gamma)}{4}\eta^{2})=0.$$
Considering Kummer's transformation
\begin{equation}
\label{relacionKummer}
M(a,b;z)=e^{z}M(b-a,b;-z)
\end{equation}
we obtain
\begin{equation}
\label{primercero}
M(\frac{N}{2}-\frac{\al}{2}-1,\frac{N}{2};\frac{(1-\gamma)}{4}\eta^{2})=0.
\end{equation}
Thus, the existence of the root $\eta_{1}$ and its behaviours are related directly with the first positive zero of
Kummer function $M(a,b;z)$ with $a=\frac{N}{2}-\frac{\al}{2}-1$. From (\ref{kummer}), we note that in cases
$\frac{N}{2}-\frac{\al}{2}-1\leq0$ the equation (\ref{primercero}) has no positive solution, obtaining  a
sufficient condition for the non-existence of $\eta_{1}$. Now, knowning that if $a < 0$ the function $M(a,N/2;z)$ has exactly $-\lfloor a\rfloor$
positive zeros,  when $\al>N-2$ the first root $\eta_{1}$ exists. $\square$

We note that the representation to $f$ given in (\ref{repre1}) is obtained when $\al\leq N-2$. To respect the inequality $\al>N-2$,  in \cite{kamin1991barenblatt} (Theorem 2.3) is used to remark that $\al_{0}>N-2$. This result is direct from  Lemma \ref{cambiodesigno1}.
\subsection{Auxiliary  results on  sign change of the profile}\label{repliminares}
In this part we present some comments and auxiliary results to understand the reasons of the oscillatory behaviors of  $f(\eta;\al)$. Moreover, the auxiliary results in this section are fundamental in our proofs to understand the representation of $f(\eta;\al)$ through the confluent hypergeometric functions and the  relation between the asymptotic representation and the similarity exponents $\al_{k}$  in (\ref{valorespropiosint}).

We continue assuming that (\ref{ecuacionparalosceros}) has exactly
 $n$  positive solutions which are denoted by
 $$\eta_{1}<\eta_{2}<\dots<\eta_{n}$$
 Let us also assume that the roots  $\eta_{m}$  can be considered as regular functions of  $\al$.
 The dependence $\eta_{m}=\eta_{m}(\al)$ is verified directly through the representations of $f_{m}(\eta;\al)$ using
 confluent hypergeometric functions, see (\ref{representacioni}) below. On the other hand, from the  representations given later in (\ref{fun2}) and (\ref{representacioni}), the function
\begin{equation}
\label{fundosvar}
h(\eta,\al)=f(\eta;\al)\qquad
\end{equation}
 is continuous for each $ (\eta,\al) \in \rr^{+}\times\rr^{+}$. This result is used to prove the following lemma.
 \begin{lemma}
 \label{cambio1}
 The functions $h_{m}(\al)=f(\eta_{m};\al)$ do not change signs.
 \end{lemma}
 \noindent{\bf{Proof:}} We begin analyzing the sign of $h_{1}(\al)$.  The proof is similar for the general case.\\
 We know that $f_{1}(\eta_{1};\al_{*})>0$ for some $\al_{*}$. If we assume $f_{1}(\eta_{1};\al^{*})<0$ for some $\al^{*}>\al_{*}$,
 through the continuity of function $h(\eta,\al)$ in (\ref{fundosvar}) there must exist some $\tilde{\al}$ in $]\al_{*},\al^{*}[$
 such that $f_{1}(\eta_{1};\tilde{\al})=0$. Directly from the definition of $\eta_{1}=\eta_{1}(\tilde{\al})$ we get $f'_{1}(\eta_{1};\tilde{\al})=0$.
 Considering
 \begin{equation*}
\left\lbrace  \begin{array}{rcl}
                f''+\frac{N-1}{\eta}f'&=&\displaystyle{-\frac{(1-\gamma)}{2} (  \tilde{\alpha} f+ \eta f  )}\\
                                 f(\eta_{1})  &=& f_{1} (\eta_{1};\tilde{\al})=0\\
                                 f'(\eta_{1}) &=& f'_{1}(\eta_{1};\tilde{\al})=0
              \end{array}
\right.
\end{equation*}
 and using existence and uniqueness results we get $f_{1}(\eta;\tilde{\al})=0$, obtaining a contradiction.
 Thus, the function $h_{1}(\al)$ does not change sign. In the general case   $m=1,2,\dots,n$, the argument follows in a similar way and therefore $h_{m}(\al)= f_{m}(\eta_{i};\al)$  do not change sign. $\square$

Knowing that $f$ is defined as piecewise function,  the next step is to study the oscillatory behavior  in each range $]\eta_{m-1},\eta_{m}[$.

\begin{lemma}
\label{signo2}
The function $f(\eta;\al)$ changes sign exactly once in $]\eta_{m-1},\eta_{m}[$.
\end{lemma}
\noindent{\bf{Proof:}} The result is obtained through contradiction method, studying the behavior of the following function
                       \begin{equation}
                       \label{funauxiliar}
                       F(\eta)=\eta f'(\eta;\al)+\al f(\eta;\al).
                       \end{equation}
                       Our analysis begins at interval
                       $]\eta_{1},\eta_{2}[$ where
                       $F(\eta)<0$. If we assume that $f(\eta;\al)$ changes sign more than once
                       in $]\eta_{1},\eta_{2}[$, then there exists at least one $\eta^{*}\in ]\eta_{1},\eta_{2}[$ such that
                       $f(\eta^{*};\al)=0$ and $f'(\eta^{*};\al)>0$. Using (\ref{funauxiliar}) we get $F(\eta^{*})>0$,
                       obtaining a contradiction.
                       Now, we assume that $f(\eta;\al)$ does not change sign in $]\eta_{1},\eta_{2}[$. Knowing that $f(\eta_{1};\al)>0$
                       we continue working under the assumption $f(\eta;\al)>0$ in $]\eta_{1},\eta_{2}[$ and therefore $f'(\eta;\al)<0$
                       in $]\eta_{1},\eta_{2}[$ (from $F(\eta)<0$ in $]\eta_{1},\eta_{2}[$).\\
                       Since (\ref{funauxiliar}), we have directly
                       \begin{equation*}
                       \begin{array}{rcl}
                       \displaystyle{\frac{dF}{d\eta}}&=&\eta f'' +(\al +1) f'\\
                                                      &=&\eta \left( f'' +\frac{N-1}{\eta}f'  \right)+(\al-(N-2))f'.
                       \end{array}
                       \end{equation*}
                       On the other hand, knowing that $f$ is a solution to (\ref{cauchynolinear}), in $]\eta_{1},\eta_{2}[$
                       we get $\sigma(\al f+\eta f')=-\frac{1+\gamma}{2}$ and therefore
                       $$\displaystyle{\frac{dF}{d\eta}=-\eta\frac{1+\gamma}{2}F +(\al-(N-2))f'}. $$
                       This equation is written as follows
                       $$ \frac{d}{d\eta}\left(  e^{\frac{(1+\gamma)}{4}\eta^{2}} F(\eta)  \right)= (\al-(N-2))e^{\frac{(1+\gamma)}{4}\eta^{2}} f'. $$
                       Integrating the equation above between $\eta_{1}$ and $\eta_{2}$ , and using $F(\eta_{1})=F(\eta_{2})=0$, we get
                        $$ \displaystyle{ \int_{\eta_{1}}^{\eta_{2}}  (\al-(N-2))e^{\frac{(1+\gamma)}{4}\eta^{2}} f' d\eta =0  }  $$
                       Under the assumption $f'(\eta;\al)<0$ in $]\eta_{1},\eta_{2}[$ and knowing $\al>N-2$, a contradiction is obtained.
                       Thus,  functions $f(\eta;\al)$ change sign exactly once in the interval  $]\eta_{1},\eta_{2}[$.\\
                       Finally, as
                       \begin{equation*}
                       \label{signoexpresion}
                       \textrm{Sign}(\eta f'_{m} +\al f_{m})=(-1)^{m+1}\qquad \textrm{when}\ \eta_{m-1}<\eta<\eta_{m},\quad m=1,2,\dots,n.
                       \end{equation*}
                       and following  previous arguments, for each $m=2,3,\dots,n$ we obtain the result in
                       $]\eta_{m-1},\eta_{m}[$. $\square$

A direct consequence of the result above  is $f(\eta_{m};\al)\cdot f(\eta_{m+1};\al)<0$. As $f(\eta_{1};\al)>0$, we get
\begin{equation}
\label{signofuncion}
\textrm{Sign}(f(\eta_{m};\al))=(-1)^{m+1}\qquad \textrm{for each}\ m=1,2,\dots,n
\end{equation}
The result in the following lemma will be used to characterize the decay rate of $f(\eta;\al)$ when $\eta\to\infty$
\begin{lemma}
\label{signo}
 In $]\eta_{n},\infty[$ the function $f(\eta;\al)$ does not change sign.
\end{lemma}
\noindent{\bf{Proof:}} The argument is similar as in the proof of Lemma \ref{signo2}. We continue using $F(\eta)$ given in (\ref{funauxiliar}) and using (\ref{signofuncion}) we know that the sign of $F(\eta)$ is determinate by the parity of $n$. We begin considering  $\al f+\eta f'<0$ in $]\eta_{n},\infty[$, i.e. $n$ odd. The proof for the even case is similar.\\
                       From (\ref{signofuncion}) we get $f(\eta_{n})>0$. If we assume that $f(\eta;\al)$ change sign
                       more than once in $]\eta_{n},\infty[$, then, there exists $\eta^{*}>\eta_{n}$ such that $f(\eta_{*};\al)=0$ and
                       $f'(\eta_{*};\al)>0$. From  the assumption  $F(\eta)<0$ a contradiction is obtained. Now, we assume that
                       $f(\eta;\al)$ changes sign exactly once in $[\eta_{n},\infty[$.\\
                       Let $\eta_{*}$ be the point where $f(\eta_{*};\al)=0$. From the proof of Lemma \ref{signo2}, $f$ can't change sign more than once, therefore $f(\eta;\al)<0$ in $]\eta_{*},\infty[$. As $f(\eta;\al)\to 0$ when $\eta\to\infty$,
                       there exists  a certain $\eta^{*}>\eta_{*}$ such that $f'(\eta^{*};\al)=0$ and $f'(\eta;\al)>0$ when $\eta>\eta^{*}$.\\
                       Knowing that $f''(\eta;\al)\to 0$ when $\eta\to \infty$, using the  ODE in (\ref{cauchynolinear}) for $\eta>\eta^{*}$ and the definition in (\ref{funauxiliar}), we get:
                       \begin{equation}
                       \label{comportamiento}
                       \lim_{\eta\to\infty}F(\eta)=0.
                       \end{equation}
                       Similar to the  proof for Lemma \ref{signo2}, considering the assumption $f'(\eta;\al)>0$ in $]\eta_{n},\infty[$
                       and knowing that $\al>N-2$, we have $(\al+1)f' > (N-1) f'$. Thus:
                       \begin{equation*}
                       \begin{array}{rcl}
                       \displaystyle{\frac{dF}{d\eta}}&>&\eta \left( f''+\frac{N-1}{\eta}f'  \right)\\
                                       &>&-\eta\frac{(1+\gamma)}{2}(\al f+\eta f).
                       \end{array}
                       \end{equation*}
                       From the assumption $F<0$ we get:
                       $$-\frac{d F}{F}>\frac{(1+\gamma)}{2}\eta d\eta.$$
                       Since $F(\eta^{*})=\al f(\eta^{*};\al)+\eta^{*}f(\eta^{*};\al)=\al f(\eta^{*};\al)<0$, integrating in
                       $[\eta^{*},\eta[$ the inequality above, we obtain:
$$ -F(\eta)>(-\al f(\eta^{*};\al))e^{-\frac{(1+\gamma)}{4}(\eta^{*})^{2}}  e^{\frac{(1+\gamma)}{4}\eta^{2}}.$$
Considering $\eta\to\infty$ we obtain a contradiction with the fact (\ref{comportamiento}).
Thus, $f(\eta;\al)$ can not change sign in $[\eta_{n},\infty[$. Finally, if we assume $\al f+\eta f'>0$ in $]\eta_{n},\infty[$,
taking $\tilde{F}(\eta)=-F(\eta)$ we obtain similar contradictions.  $\square$

\section{Representation of the similarity solutions}\label{sectionsimilarity}
In this part we detail the representations of $f$ using the Kummer and Tricomi functions.\\
For simplicity, the method to describe the profiles $f(\eta;\al)$ is developed by separate according to the following
cases: positive  and sign change similarity solutions.
\subsection{Positive similarity solutions}
We continue working under the assumption $\al>N-2$, therefore $\eta_{1}$ exists (see Lemma \ref{cambiodesigno1}), in other case  the representation of $f$ is given by (\ref{repre1}).

We begin considering the second range given in (\ref{cauchylinear2}). Taking $m=2$ in (\ref{transformacion}), $f_{2}(\eta;\al)$ can be represented by two linearly independent solutions of Kummer equations.\\
Knowing that (\ref{kummer}) is the first solution, in this part we consider the Tricomi function $U(\al/2,N/2;z)$ as the second linearly independent solution of (\ref{kummereq}).\\
As (see formulae  13.1.22 \cite{abramowitz1972handbook})
\begin{equation}
\label{wronskianogen}
\mathcal{W}\lbrace M(\al/2,N/2;z)),U(\al/2,N/2;z) \rbrace =-\frac{\Gamma(N/2)}{\Gamma(\al/2)} z^{-N/2}e^{z},
\end{equation}
taking $\varphi_{1}(\eta)=M(\al/2,N/2;-\frac{(1+\gamma)}{4}\eta^{2})$ and $\varphi_{2}(\eta)=U(\al/2,N/2;-\frac{(1+\gamma)}{4}\eta^{2})$,
complex solutions for the ODE in (\ref{cauchylinear2}),  the solution $f_{2}(\eta;\al)$ can be written as follows:
\begin{equation}
\label{representacion2}
f_{2}(\eta;\al)=A_{2}(\al) M(\frac{\al}{2},\frac{N}{2};-\frac{(1+\gamma)}{4}\eta^{2}) + B_{2}(\al) U(\frac{\al}{2},\frac{N}{2};-\frac{(1+\gamma)}{4}\eta^{2}).
\end{equation}
Here, we consider the principal branch for $U(\al/2,N/2;z)$, i.e. $ph z\in[0,\pi]$. Knowing that $f_{2}(\eta;\al)$  must satisfy $f_{2}(\eta_{1};\al)=f_{1}(\eta_{1};\al)$ and
 $f_{2}'(\eta_{1};\al)=f_{1}'(\eta_{1};\al)$, to determinate $A_{2}(\al)$ and $B_{2}(\al)$ we consider
 \begin{equation*}
 \begin{pmatrix} M(\frac{\al}{2},\frac{N}{2};s_{1}) & U(\frac{\al}{2},\frac{N}{2};s_{1}) \\
                 M'(\frac{\al}{2},\frac{N}{2};s_{1})& U'(\frac{\al}{2},\frac{N}{2};s_{1}) \end{pmatrix}
                 \begin{pmatrix} A_{2}(\al)\\ B_{2}(\al)  \end{pmatrix}= \begin{pmatrix} f_{1}(\eta_{1};\al)\\ -\frac{2}{(1+\gamma)\eta_{1}} f_{1}'(\eta_{1};\al) \end{pmatrix},
 \end{equation*}
 with $s_{1}=-\frac{(1+\gamma)}{4}\eta^{2}$. On the other hands, as $\al f_{1}(\eta_{1};\al)+\eta_{1}f_{1}'(\eta_{1}\al)=0$, hence
\begin{equation*}
A_{2}(\al)=\frac{f_{1}(\eta_{1};\al )}{s_{1}\mathcal{W}_{1}}( s_{1} U'(\frac{\al}{2},\frac{N}{2};s_{1}) +\frac{\al}{2} U(\frac{\al}{2},\frac{N}{2};s_{1}))
\end{equation*}
and
\begin{equation*}
B_{2}(\al)=\frac{f_{1}(\eta_{1};\al )}{s_{1}\mathcal{W}_{1}}( s_{1} M'(\frac{\al}{2},\frac{N}{2};s_{1}) +\frac{\al}{2} M(\frac{\al}{2},\frac{N}{2};s_{1})),
\end{equation*}
where  $\mathcal{W}_{1}$  denote the Wronskian (\ref{wronskianogen}) evaluated at $\eta=\eta_{1}$.\\
The next step is to characterize the constants $A_{2}(\al)$ and $B_{2}(\al)$. We begin with $A_{2}(\al)$, where  are used recurrence results of $U(a,b;z)$.\\ From the relation (see  \cite{abramowitz1972handbook} 13.4.25)
$$ aU(a,b;z)+zU'(a,b;z)=a(1+a-b)U(a+1,b;z), $$
follows
\begin{equation}
\label{constanteA}
A_{2}(\al)=\frac{f_{1}(\eta_{1};\al )}{s_{1}\mathcal{W}_{1}}\frac{\al}{2}\left\lbrace 1+\frac{\al-N}{2} \right\rbrace U(\frac{\al}{2}+1,\frac{N}{2};s_{1}).
\end{equation}
Now, we use recurrence relations of the Kummer function to characterize $B_{2}(\al)$.
From (\ref{recurrenciaKummer}), we get
\begin{equation}
\label{constanteB}
B_{2}(\al)=-\frac{f_{1}(\eta_{1};\al )}{s_{1}\mathcal{W}_{1}}\frac{\al}{2}M(\frac{\al}{2}+1,\frac{N}{2};s_{1}).
\end{equation}
Therefore, considering $z=-\frac{(1-\gamma)}{4}\eta^{2}$, while $f>0$, its representation is given by
\begin{equation}
\label{fun2}
f(\eta;\al)=\left\lbrace  \begin{array}{rl}
                           M(\frac{\al}{2},\frac{N}{2}; z)&  0\leq \eta\leq \eta_{1}\\
                           A_{2}(\al) M(\frac{\al}{2},\frac{N}{2};\frac{1+\gamma}{1-\gamma} z) + B_{2}(\al) U(\frac{\al}{2},\frac{N}{2};\frac{1+\gamma}{1-\gamma} z)& \eta_{1}\leq \eta
                          \end{array}
            \right.
\end{equation}
where $A_{2}(\al)$ and $B_{2}(\al)$ are defined in (\ref{constanteA}) and (\ref{constanteB}) respectively.

We note that,  if  $\eta_{1}$ is the unique solution to (\ref{ecuacionparalosceros}), we have $f>0$. In other case we consider the result presented in Lemma \ref{signo2}.  This is proved directly  by the asymptotic representations of $f(\eta;\al)$ shown in the next section.

\subsection{Similarity solutions with sign changes}
To illustrate our method, we present the case $n=2$ in (\ref{formak}),thus,  $\al f+\eta f'=0$ has two roots  denoted by $\eta_{1}<\eta_{2}$.  Following similar ideas, the remaining cases are presented below.\\
Let $f(\eta;\al)$ be the solution of problem (\ref{cauchynolinear}) given in (\ref{formak}) with $n=2$.\\
Using the representation of $f_{2}(\eta;\al)$  in (\ref{representacion2}), Eq. (\ref{ecuacionparalosceros})
is written as follows
\begin{equation}
\label{raiz2}
\begin{array}{l}
\al\left\lbrace A_{2}(\al)M(\frac{\al}{2},\frac{N}{2};s_{2}) + B_{2}(\al)U(\frac{\al}{2},\frac{N}{2};s_{2})  \right\rbrace- \\
 \\
 -s_{2}\left\lbrace  A_{2}(\al)M'(\frac{\al}{2},\frac{N}{2};s_{2}) + B_{2}(\al)U'(\frac{\al}{2},\frac{N}{2};s_{2})  \right\rbrace=0,
\end{array}
\end{equation}
with $s_{2}=-\frac{(1+\gamma)}{4}(\eta_{2})^{2}$. Following similar ideas as those shown in the descriptions of the constants $A_{1}(\al)$ and $B_{1}(\al)$,
we prove that  equation (\ref{raiz2}) can be written
\begin{equation}
\label{raiz2b}
A_{2}(\al) M(\frac{\al}{2}+1,\frac{N}{2};s_{2}) = - B_{2}(\al)\left\lbrace 1+\frac{\al-N}{2}  \right\rbrace U(\frac{\al}{2}+1,\frac{N}{2};s_{2}).
\end{equation}
Considering $\eta_{2}$, the second root of $\eta f'_{2}(\eta;\al)+\frac{\al}{2}f_{2}(\eta;\al)=0$, we get $s_{2}=-\frac{(1-\gamma)}{4}(\eta_{2})^{2}$ the root in (\ref{raiz2b}). Notice that $s_{2}$ is described using the first root $\eta_{1}$ and $\al$.  From (\ref{primercero}) and using the implicit function theorem,  $\eta_{1}$ can be defined as regular function $\eta_{1}=\eta_{1}(\al)$ and therefore $\eta_{2}=\eta_{2}(\al)$.

We know that $f(\eta;\al)=f_{3}(\eta;\al)$ in $[\eta_{2},\infty[$. Using (\ref{transformacion}) follows
\begin{equation}
\label{representacion3}
f_{3}(\eta;\al)=A_{3}(\al) M(\frac{\al}{2},\frac{N}{2};-\frac{(1-\gamma)}{4}\eta^{2}) + B_{3}(\al) U(\frac{\al}{2},\frac{N}{2};-\frac{(1-\gamma)}{4}\eta^{2})
\end{equation}
Hence, following similar arguments as in the analysis for $f_{2}$, we get:
\begin{equation}
\label{constantes3}
\begin{array}{rcl}
A_{3}(\al)&=&\displaystyle{\frac{f_{2}(\eta_{2};\al )}{s_{2}\mathcal{W}_{2}}\frac{\al}{2}\left\lbrace 1+\frac{\al-N}{2} \right\rbrace U(\frac{\al}{2}+1,\frac{N}{2};s_{2})}\\
B_{3}(\al)&=&\displaystyle{-\frac{f_{2}(\eta_{2};\al )}{s_{2}\mathcal{W}_{2}}\frac{\al}{2}M(\frac{\al}{2}+1,\frac{N}{2};s_{2})},
\end{array}
\end{equation}
where $\varphi_{1}(\eta)=M(\al/2,N/2;-\frac{(1-\gamma)}{4}\eta^{2})$, $\varphi_{2}(\eta)=U(\al/2,N/2;-\frac{(1-\gamma)}{4}\eta^{2})$, and $\mathcal{W}_{2} $ is the Wronskian between $\varphi_{1}$ and $\varphi_{2}$ at $\eta=\eta_{2}$.\\
On the other hand, since $A_{3}(\al)$ depends of $\eta_{2}$, then $A_{3}(\al)$ depends of $\eta_{1}$ (similar to $B_{3}(\al)$).
For the general case, consider the $n$ positive roots:
$\eta_{1}<\eta_{2}<\dots<\eta_{n}$
of  equation (\ref{ecuacionparalosceros}). For each $m=1,2,\dots, n+1$, consider $f_{m}(\eta;\al)$ as (\ref{formak}).  $f_{1}(\eta;\al)$ is defined by (\ref{repre1}) and in the cases $m\geq 2$, we have
\begin{equation}
\label{representacioni}
f_{m}(\eta;\al)=A_{m}(\al) M(\frac{\al}{2},\frac{N}{2};z_{m})+ B_{m}(\al) U(\frac{\al}{2},\frac{N}{2};z_{m}),
\end{equation}
where $z_{m}=-\frac{(1+(-1)^{m}\gamma)}{4}\eta^{2}$ with $ \eta_{m-1}\leq \eta\leq\eta_{m}$. Likewise, we define  $f_{n+1}(\eta;\al)$ for $\eta\in[ \eta_{n},\infty[$. In a similar way as in the statement of (\ref{constantes3}), for the general case, we have:
\begin{equation}
\label{constantesi}
\begin{array}{rcl}
A_{m}(\al)&=&\displaystyle{\frac{f_{m-1}(\eta_{m-1};\al)}{s_{m-1}\mathcal{W}_{m-1}}\frac{\al}{2}\left\lbrace 1+\frac{\al-N}{2} \right\rbrace U(\frac{\al}{2}+1,\frac{N}{2};s_{m-1})}\\
\\
B_{m}(\al)&=&-\displaystyle{\frac{f_{m-1}(\eta_{2};\al)}{s_{m-1}\mathcal{W}_{m-1}}\frac{\al}{2}M(\frac{\al}{2}+1,\frac{N}{2};s_{m-1})},
\end{array}
\end{equation}

where $s_{m-1}=-\frac{(1+(-1)^{m})}{4}(\eta_{m-1})^{2}$ and $\mathcal{W}_{m-1}$ the Wronskian (\ref{wronskianogen})  at $\eta_{m-1}$. Here, $\eta_{m-1}$ is the $m-1$ root of $\eta f_{m-1}+\frac{\al}{2}f_{m-1}=0$. These root is determined from
\begin{equation*}
\label{raizn}
A_{m-1}(\al) M(\frac{\al}{2}+1,\frac{N}{2};s_{m-1}) = - B_{m-1}(\al)\left\lbrace 1+\frac{\al-N}{2}  \right\rbrace U(\frac{\al}{2}+1,\frac{N}{2};s_{m-1}).
\end{equation*}
 We remark that through the implicit function theorem, the dependence of $A_{m-1}$ with respect to similarity parameter $\al$ is obtained. This dependence is obtained by recurrence, using the results in each previous
interval. That is, $A_{m-1}$ depends of $A_{k}$; $B_{k}$; $s_{k}$ with $k=1,2,\dots,m-2$. Similar for $B_{m-1}$ and $s_{m-1}$.

The  asymptotic behaviors of $f(\eta;\al)$ are obtained using the asymptotic representations
of confluent hypergeometric functions $M(a,b;z)$ and $U(a,b;z)$. These ideas are developed in the next section.


\section{Main Results} \label{principal}
Since our profiles $f(\eta;\al)$ can be defined from (\ref{formak}), their asymptotic behaviors are described by the asymptotic behavior of the Kummer and Tricomi functions, see (\ref{representacioni}).\\
We begin considering the asymptotic representation of the Kummer function given by
\begin{equation}
\label{kummerasym}
\begin{array}{lr}
M(a,b;z)=\displaystyle{\frac{\Gamma(b)}{\Gamma(b-a)}(-z)^{-a}\left[ \sum_{k=0}^{N}\frac{(a)_{k}(a-b+1)_{k}}{k! (-z)^{k}} +O(z^{-N}) \right]}\quad (z<0).
\end{array}
\end{equation}
This asymptotic result is valid when $b-a\neq -l$ with $l=0,1,2,\dots$, while in cases $b-a=-l$ the behavior is obtained from (\ref{decaimiento}).\\
Throughout our proofs, we consider the following confluent hypergeometric functions (see \cite{olver1974introduction} 10.09)
\begin{equation}
\label{kummer2}
\begin{array}{rcl}
V(a,b;z) &=&\displaystyle{\frac{\Gamma(a)}{\Gamma(b)}e^{(b-a)i\pi}M(a,b;z) -\frac{\Gamma(a)}{\Gamma(b-a)} e^{bi \pi }U(a,b;z)}.
\end{array}
\end{equation}

The asymptotic representations of our  profiles $f(\eta;\al)$ are also written in terms of $V(a,b;z)$. We note that the asymptotic representation for $V(a,b;z)$ is given by (see  10.02 in \cite{olver1974introduction} )
\begin{equation}
\label{tricomiasym}
V(a,b;z)= e^{z}(-z)^{a-b}\displaystyle{  \left[ \sum_{k=0}^{N}\frac{(b-a)_{k} (1-a)_{k}}{k! z^{k}} +O(z^{-N}) \right]}\qquad (z<0).
\end{equation}
Therefore,  using (\ref{kummerasym}) and (\ref{tricomiasym}), the asymptotic representation for $U(a,b;z)$ can be obtained from (\ref{kummer2}).\\
On the other hand, the function in (\ref{kummer2}) is also introduced as follows (see 10.03 in \cite{olver1974introduction}):
\begin{equation}
\label{kummer3}
V(a,b;z)=e^{z}U(b-a,b;z)
\end{equation}
From this definition, we note that the gaussian-type decays of $f(\eta;\al)$ are related with $V(a,b;z)$.

\noindent{\bf Proof Theorem \ref{principal1}:}  To illustrate our method, we begin studying the existence of $\al=\al_{0}$, i.e. $\eta_{1}$ is the unique solution for the equation $\eta f'+\al f=0$. The general case is studied in a similar form. In this proof we assume that $\al\neq N$, the other case is analyzed in a similar way, but considering  other two linearly independent solution for the Kummer equation.\\
When $\eta_{1}$ is the unique solution of $\eta f'+\al f=0$ and knowing that $f(\eta;\al)=f_{2}(\eta;\al)$ for $\eta>\eta_{1}$,  follows that the asymptotic behavior of $f(\eta;\al)$ is given by the asymptotic representation of $f_{2}$.\\
Assuming $\displaystyle{\Gamma(N/2)A_{2}(\al)+e^{-i\pi \al/2}\Gamma(\lbrace N-\al  \rbrace/2) B_{2}(\al)\neq0}$, from the asymptotic behaviors given in (\ref{kummerasym}) and (\ref{tricomiasym})  when $\eta$ is large we get the following asymptotic representation:
\begin{equation}
\label{asintotico1}
f_{2}(\eta;\al)\sim  \left\lbrace \frac{\Gamma(\frac{N}{2})}{\Gamma(\frac{N}{2}-\frac{\al}{2} )}A_{2}(\al)+e^{-i\frac{\pi}{2} \al} B_{2}(\al)  \right\rbrace  \left(\frac{1+\gamma}{4} \eta^{2}\right)^{-\al/2}.
\end{equation}
If $\eta_{1}$ is the unique root, then we continue  analyzing the asymptotic representation for $f_{2}(\eta;\al)$ (the general case is studied below).\\
Using  (\ref{constanteA}) and (\ref{constanteB}), the asymptotic representation (\ref{asintotico1}) is written as follows:
\begin{equation*}
\begin{array}{ll}
f_{2}(\eta;\al)\sim& -\displaystyle{ \frac{f_{1}(\eta_{1};\al)}{s_{1}\mathcal{W}_{1}}\frac{\al}{2}\left\lbrace  \frac{\Gamma(\frac{N}{2})}{\Gamma(\frac{N}{2}-\frac{\al}{2})} \left( 1+\frac{\al}{2}-\frac{N}{2} \right)U(\frac{\al}{2}+1,\frac{N}{2};s_{1})  \right.}\\
\\
&\displaystyle{\left.+ e^{-i\frac{\pi}{2} \al}  M(\frac{\al}{2}+1,\frac{N}{2};s_{1})  \right\rbrace \left((1+\gamma)\frac{\eta^{2}}{4}\right)^{-\frac{\al}{2}}}.
\end{array}
\end{equation*}
Now, from the definition in (\ref{wronskianogen}) follows
\begin{equation*}
\begin{array}{rl}
f_{2}(\eta;\al)&\sim\displaystyle{ -e^{-s_{1}}(s_{1})^{\frac{N}{2}-1}f_{1}(\eta_{1};\al) \frac{\al}{2} \left\lbrace \frac{\Gamma(\frac{\al}{2})(1+\frac{\al}{2}-\frac{N}{2})}{\Gamma(\frac{N}{2}-\frac{\al}{2})} U(\frac{\al}{2}+1,\frac{N}{2};s_{1})\right.}\\
\\
   &\displaystyle{\left.-\frac{e^{-i\frac{\pi}{2}\al }\Gamma(\frac{\al}{2})}{\Gamma(\frac{N}{2})} M(\frac{\al}{2}+1,\frac{N}{2};s_{1})     \right\rbrace \left((1+\gamma)\frac{\eta^{2}}{4}\right)^{-\frac{\al}{2}}},
\end{array}
\end{equation*}
with $s_{1}=-(1+\gamma)\eta^{2}/4$.\\
Knowing $\frac{\al}{2}\Gamma(\frac{\al}{2})=\Gamma(1+\frac{\al}{2})$ and $\frac{1+\frac{\al-N}{2}}{ \Gamma(\lbrace  N-\al\rbrace/2)}=-\frac{1}{\Gamma(\frac{N}{2}-\frac{\al}{2}-1)}$, we get
\begin{equation*}
\begin{array}{rl}
f_{2}(\eta;\al)\sim & \displaystyle{e^{\pi i}(s_{1})^{\frac{N}{2}-1} f_{1}(\eta_{1};\al)e^{-i\frac{N}{2}\pi} \left\lbrace \frac{\Gamma(1+\frac{\al}{2})e^{i\frac{N}{2}\pi}}{\Gamma(\frac{N}{2}-\frac{\al}{2}-1)} U(\frac{\al}{2}+1,\frac{N}{2};s_{1})+\right.}\\
\\
   &\displaystyle{\left.\frac{e^{i\lbrace\frac{N}{2}-\frac{\al}{2}-1\rbrace\pi }\Gamma(\frac{\al}{2})}{\Gamma(\frac{N}{2})} M(\frac{\al}{2}+1,\frac{N}{2};s_{1})     \right\rbrace e^{-s_{1}} \left((1+\gamma)\frac{\eta^{2}}{4}\right)^{-\frac{\al}{2}}}.
\end{array}
\end{equation*}
Finally, from (\ref{kummer2}) and (\ref{kummer3}) we obtain
\begin{equation}
\label{asintotico2}
f_{2}(\eta;\al)\sim  C_{2}(\al) \eta^{-\al},\qquad\textrm{as}\ \eta\to\infty,
\end{equation}
where
\begin{equation*}
 C_{2}(\al)=((1+\gamma)/4)^{\frac{N}{2}-\frac{\al}{2}-1}\eta_{1}^{N-2} f_{1}(\eta_{1};\al) U(\frac{N-\al}{2}-1,\frac{N}{2};\frac{(1+\gamma)}{4}\eta^{2}_{1}).
\end{equation*}
Now, we study the existence of similarity exponent $\al_{0}$ through the sign change of $f_{2}(\eta;\al)$. From Lemma \ref{cambio1},  for each $\al>0$ we have $f(\eta_{1};\al)>0$ and therefore $f_{1}(\eta_{1};\al)>0$. Thus, the change in the asymptotic behavior (\ref{asintotico2}) is given by the condition
\begin{equation}
\label{ecuacion}
U(\frac{N-\al}{2}-1,\frac{N}{2};\frac{(1+\gamma)}{4}\eta^{2}_{1})=0
\end{equation}
On the other hand, we know that while $\eta_{1}$ be the unique root of $\eta f'+\al f=0$, we have $f(\eta;\al)=f_{2}(\eta;\al)$ for $\eta>\eta_{1}$. Thus, the change of asymptotic behavior of $f(\eta;\al)$ is related with (\ref{ecuacion}). Moreover, the existence of $\eta_{2}$, the second root of $\eta f'+\al f=0$, is also related with (\ref{ecuacion}), see Lemma \ref{signo}.\\
Let $\al^{*}$ be a root of (\ref{ecuacion}) with $\eta_{1}=\eta_{1}(\al)$ from (\ref{primercero}). From (\ref{asintotico1}), we have $ \frac{\Gamma(N/2)}{\Gamma(\lbrace N-\al \rbrace/2)}A_{2}(\al^{*})+e^{-i\pi \al^{*}/2} B_{2}(\al^{*})=0 $ and therefore
\begin{equation*}
A_{2}(\al^{*})=-e^{-i\pi \al^{*}/2} \frac{\Gamma(\lbrace N-\al \rbrace/2)}{\Gamma(N/2)}B_{2}(\al^{*}),
\end{equation*}
Likewise the deduction of (\ref{asintotico2}), we get
\begin{equation*}
\label{represetacion2}
f_{2}(\eta;\al^{*})=-\displaystyle{B_{2}(\al^{*})\frac{\Gamma(\lbrace N-\al^{*} \rbrace/2)}{\Gamma(\al^{*}/2)}e^{-i\frac{N}{2}\pi}V(\frac{\al}{2},\frac{N}{2};-\frac{1+\gamma}{4}\eta^{2})}
\end{equation*}
and, knowing the asymptotic representation (\ref{tricomiasym}), we obtain
\begin{equation*}
\label{asinoticorama2}
f_{2}(\eta;\alp)=D_{2}(\alp) \eta^{\alp-N}e^{-\frac{1+\gamma}{4}\eta^{2}}\left[ 1+O(\eta^{-2})  \right],\qquad\eta\to\infty,
\end{equation*}
where
\begin{equation*}
D_{2}(\alp)=-(\lbrace 1+\gamma \rbrace/4)^{\frac{\alp-N}{2}} \frac{\Gamma(\lbrace N-\alp \rbrace/2)}{\Gamma(\alp/2)} e^{-i\frac{N}{2}\pi} B_{2}(\alp)
\end{equation*}
Finally, since $\al^{*}$ is a root of (\ref{ecuacion}) and knowing the representation of $B_{2}(\al)$ in (\ref{constanteB}), we obtain $f_{2}(\eta;\alp)>0$ considering the smallest root such that
\begin{equation}
 \label{signoexponencial}
 M(\frac{N-\alp}{2}-1,\frac{N}{2};\frac{(1+\gamma)}{4}\eta_{1}^{2}) \cdot \Gamma(\lbrace N-\alp \rbrace/2)>0
\end{equation}
Under the condition (\ref{signoexponencial}), for each $0\leq \al\leq \alp$ we have $f(\eta;\al)=f_{2}(\eta;\al)$ when $\eta\geq \eta_{1}$. Taking $\al_{0}=\al^{*}$, our first critical exponent $\al=\al_{0}$ is obtained. Hence, if $\al<\al_{0}$ the asymptotic representation is given in (\ref{asintotico2}), while the asymptotic behavior at $\al=\al_{0}$ is (\ref{asinoticorama2}). For $\al=N$, from (\ref{decaimiento}),  $M(\frac{\al}{2},\frac{N}{2};-\frac{(1+\gamma)}{4}\eta^{2})$ has gaussian-type behavior  and therefore the asymptotic representation at $f_{2}(\eta;\al)$ is given in (\ref{asintotico1}) considering $A_{2}(\al)=0$. The comment above is necessary for the description of $f(\eta;\al)$ when $\al_{0}>N$, i.e. when  $\gamma>0$.

The general case is studied in a similar way. For the existence of $\al_{n}$ we assume that $\eta f'+\al f=0$ has exactly $n+1$ roots $\eta_{1}<\eta_{2}<\dots<\eta_{n+1}$.   Through the representation of $f(\eta;\al)$ given in (\ref{formak}), we note that the existence of $\al_{n}$ is obtained following a similar argument like case $\al_{0}$.  Under the assumption:
$$  \frac{\Gamma(N/2)}{\Gamma(\lbrace N-\al  \rbrace/2)}A_{n+2}(\al)+e^{-i\pi \al/2} B_{n+2}(\al)\neq 0, $$
and considering that $N-\al\neq -l$ $(l=0,1,2,\dots)$, the asymptotic representation of $f_{n+2}(\eta;\al)$ is given by
\begin{equation*}
\begin{array}{l}
f_{n+2}(\eta;\al)\sim \left(\frac{1+(-1)^{n}\gamma}{4}\right)^{-\frac{\al}{2}} \left\lbrace \frac{\Gamma(\frac{N}{2})}{\Gamma(\frac{N}{2}-\frac{\al}{2})}A_{n+2}(\al)+e^{-i\frac{\pi}{2} \al} B_{n+2}(\al)  \right\rbrace  \eta^{-\al}
\end{array}
\end{equation*}
Assuming $N-\al= -l$ $(l=0,1,2,\dots)$, we consider $A_{n+2}(\al)=0$. From the representations given in (\ref{constantesi})
 and (\ref{kummer2}), we obtain
\begin{equation*}
\label{asintoticom2}
f_{n+2}(\eta;\al)\sim  C_{n+2}(\al) \left(\frac{1+\gamma}{4}\right)^{\frac{N}{2}-\frac{\al}{2}-1}\eta^{-\al},\qquad\textrm{as}\ \eta\to\infty
\end{equation*}
where
\begin{equation*}
 C_{n+2}(\al)=\eta_{n+1}^{N-2} f_{n+1}(\eta_{n+1};\al) U(\frac{N}{2}-\frac{\al}{2}-1,\frac{N}{2};\frac{(1+(-1)^{n}\gamma)}{4}\eta^{2}_{n+1}).
\end{equation*}
In case
$$ \frac{\Gamma(N/2)}{\Gamma(\lbrace N-\al  \rbrace/2)}A_{m+2}(\al)+e^{-i\pi \al/2} B_{m+2}(\al)=0 $$
following similar step as the deduction of (\ref{represetacion2}), we get
\begin{equation*}
\label{represetacionn2}
f_{n+2}(\eta;\al^{*})=-\displaystyle{B_{n+2}(\al^{*})\frac{\Gamma(\lbrace N-\al^{*} \rbrace/2)}{\Gamma(\al^{*}/2)}e^{-i\frac{N}{2}\pi}V(\frac{\al}{2},\frac{N}{2};-\frac{1+(-1)^{n}\gamma}{4}\eta^{2})}
\end{equation*}
and, knowing the asymptotic representation (\ref{tricomiasym}), we obtain
\begin{equation*}
\label{asinoticoraman2}
f_{n+2}(\eta;\alp)=D_{n+2}(\alp) \eta^{\alp-N}e^{-\frac{1+(-1)^{n}\gamma}{4}\eta^{2}}\left[ 1+O(\eta^{-2})  \right],\qquad\eta\to\infty
\end{equation*}
where
\begin{equation*}
D_{n+2}(\alp)=-(\lbrace 1+\gamma \rbrace/4)^{\frac{\alp-N}{2}} \frac{\Gamma(\lbrace N-\alp \rbrace/2)}{\Gamma(\alp/2)} e^{-i\frac{N}{2}\pi} B_{n+2}(\alp)
\end{equation*}
Hence, knowing the representations given in (\ref{constantesi}), and assuming that $\al_{n-1}$ exist, the critical exponent $\al_{n}$  is obtained considering the smallest $\al^{*}$ such that
\begin{equation}
\label{sistema}
\begin{array}{rcl}
M(\frac{N-\alp}{2}-1,\frac{N}{2};\frac{(1-\gamma)}{4}\eta_{1}^{2})&=&0\\
U(\frac{N-\alp}{2}-1,\frac{N}{2};\frac{(1+(-1)^{n}\gamma)}{4}\eta^{2}_{n+1})&=&0\\
M(\frac{N-\alp}{2}-1,\frac{N}{2};\frac{(1+(-1)^{n}\gamma)}{4}\eta_{n+1}^{2}) \cdot \Gamma(\lbrace N-\alp \rbrace/2)\cdot (-1)^{n}&>&0,
\end{array}
\end{equation}
satisfying $\al_{n-1}<\al^{*}$. Finally, if $\al$ satisfy $\al_{n-1}<\al< \al_{n}$, the asymptotic representation is given by
(\ref{asintoticom2}), finishing our proof.   $\square$

\noindent{\bf Proof Theorem \ref{principal2}:} Consider $\al$ such that $\al_{k}<\al\leq \al_{k+1}$. Using the analysis above,
$\eta f'+\al f=0$ has exactly $k+1$ roots $\eta_{1}<\eta_{2}<\dots<\eta_{k+1}$. Thus, the profile $f(\eta;\al)$ is described through $k+2$ linear Cauchy problems, each  defined in $[0,\eta_{1}[;\  [\eta_{1},\eta_{2}[;\dots[\eta_{k},\eta_{k+1}[;\ [\eta_{k+1},\infty[$, respectively.
From Lemma \ref{signo2},  $f(\eta;\al)$ has exactly $k$ zeros in $[0,\eta_{k+1}[$. Finally, using Lemma \ref{signo} the profile $f(\eta;\al)$ does not change sign in $[\eta_{n+1},\infty[$,
and therefore the proof is finished. $\square$


\section{Extension to a nonlinear diffusion equation}\label{extension}
In this section we show how the results presented in Theorems \ref{principal1} and \ref{principal2} can be extended
on certain parabolic Bellman equations
\begin{equation}
\label{unif}
u_{t}=F(D^{2}u)\qquad \textrm{in}\ \rn\times]0,\infty[.
\end{equation}
Here, the diffusion term is defined by a mapping $F:\sm\to\rr$, where $\sm$  denotes the space of symmetric matrices
of order $N\times N$ and $D^{2}u$  denotes the Hessian matrix of $u(x,t)$ with respect to the spacial variable $x\in\rn$.\\
In regard of the general problem (\ref{unif}), a complete analysis on asymptotic behaviors for $t$ large of the solutions  was  presented in \cite{armstrong2010long}, where  relation between asymptotic behaviors and self-similar solutions of (\ref{unif})
with gaussian-type decays was obtained. The conditions to this asymptotic behavior is given by the assumption that the positive initial conditions $u(x,0)$ has spatial decay at most of gaussian-type.

Analysis on (\ref{unif}) can begin considering the Pucci extremal operators. These extremal operators can be introduced as follows (see \cite{luis1995fully}):
\begin{equation}
\label{extremalpucci}
\mathcal{M}_{\lambda,\Lambda}^{-}(A)=\la \tr(A^{+})-\Lambda\tr (A^{-}),\qquad \mathcal{M}_{\lambda,\Lambda}^{+}(A)=\Lambda \tr(A^{+})-\la\tr (A^{-}),
\end{equation}
with $A^{-},\ A^{+}$ semi-defined positives matrices, such that $A=A^{+}-A^{-}$, $A^{-}\cdot A^{+}=0$, i.e., the orthogonal
decomposition of $A$.\\
Thus,  when $F$ is uniformly elliptic and positively homogeneous operator, studying the following uniformly parabolic equation
\begin{equation}
\label{unif2}
u_{t}=\mathcal{M}_{\lambda,\Lambda}^{+}(D^{2}u)\qquad \textrm{in}\ \rn\times]0,\infty[,
\end{equation}
we obtain super-solutions for (\ref{unif}). Similarly, sub-solutions are obtained considering $\mathcal{M}_{\lambda,\Lambda}^{-}$. Hence, using comparison results, many qualitative behaviors for (\ref{unif}) can be studied through the solutions of (\ref{unif2}).\\
If we consider that (\ref{unif2})
admits solutions of the form (\ref{selfsimilar}), the profiles $f(\eta;\al)$ satisfy the following nonlinear Cauchy problem:
\begin{equation}
\label{cauchynolinear2}
\left\lbrace  \begin{array}{rcl}
                \sigma(f'')f''+\sigma(f')\frac{N-1}{\eta}f'&=& \frac{\alpha}{2} f'+ \frac{\eta}{2} f  \\
                                 f(0) &=&1\\
                                 f'(0)&=&0
              \end{array}
\right.
\end{equation}
Here $\sigma(z)$ is a discontinuous function defined as follows: $\sigma(z)=\Lambda$ if $z>0$ and $\sigma(z)=\lambda$
if $z<0$. The analysis of (\ref{cauchynolinear2}) for the case $\al=\al_{0}$
was studied in \cite{kroger1994regularity}.\\
The problem (\ref{cauchynolinear2}) is analyzed by linear problems, similar to the problem (\ref{cauchynolinear}).\\
For (\ref{cauchynolinear2}), the ranges of each linear problem are defined using roots $\tilde{\eta}_{k}$
of $f''(\eta)=0$ and the roots $\hat{\eta}_{k}$ of   $f'(\eta)=0$.
Being careful with some details, the solutions of (\ref{cauchynolinear2}) can be determined and represented by confluent
hypergeometric functions, having powerful tools to obtain results similar as to  Theorems \ref{principal1}
and \ref{principal2}.
For $\mathcal{M}_{\lambda,\Lambda}^{-}$ the result is similar considering in (\ref{cauchynolinear2}) the function
$\sigma$ such that $\sigma(z)=\lambda$ if $z>0$ and $\sigma(z)=\Lambda$ if $z<0$.

\section{Comments and Conclusions}
In this article we developed a method to describe the similarity solution to (\ref{barenblattcanonico}). Using a simple change of variable we obtain an explicit representation of the similarity solution through the confluent hypergeometric functions, obtaining a generalization of the results presented in \cite{barenblatt1969self,Kerchman1971self}. We applied the explicit representation  to describe the oscillatory and asymptotic behaviors of the similarity solution, obtaining similar results to the convection heat equation (case $\gamma=0$). The gaussian-type decays  show the generation of a new zero
for the profiles. This feature can be seen in asymptotic representation (\ref{asintoticom2}). Each  exponent $\al_{0}<\al_{1}<\dots$ in (\ref{valorespropiosint}) is related with the zeros of Kummer and Tricomi functions. This analysis allows generating new explicit approximations for anomalous exponents,
which we hope to develop in a future paper.\\\
Our method can be easy applied on other nonlinear parabolic equations.  We believe that through the asymptotic behaviors of similarity solution of (\ref{unif2}) it's possible  to obtain  conditions on $u_{0}$ for the existence of global solution to fully nonlinear problems
\begin{equation*}
\begin{array}{rcll}
            u_{t} &=&F(D^{2}u)+u^{p}&x\in\rn,\ t>0,\ u\geq 0,\ p>1\\
            u(x,0)&=&u_{0}(x)& u_{0}\geq 0,\ u_{0}\neq0
\end{array}
\end{equation*}
These topics are beyond to scope of this article and we leave them for future studies.

\bibliographystyle{plain}
\bibliography{biblio}{}

\begin{thebibliography}{10}

\bibitem{abramowitz1972handbook}
Milton Abramowitz and Irene~A Stegun.
\newblock {\em Handbook of mathematical functions: with formulas, graphs, and
  mathematical tables}.
\newblock Number~55. Courier Dover Publications, 1972.

\bibitem{armstrong2010long}
Scott~N Armstrong and Maxim Trokhimtchouk.
\newblock Long-time asymptotics for fully nonlinear homogeneous parabolic
  equations.
\newblock {\em Calculus of Variations and Partial Differential Equations},
  38(3-4):521--540, 2010.

\bibitem{aronson1994calculation}
DG~Aronson and JL~Vazquez.
\newblock Calculation of anomalous exponents in nonlinear diffusion.
\newblock {\em Physical review letters}, 72(3):348, 1994.

\bibitem{barenblatt1955concerning}
GI~Barenblatt and AP~Krylov.
\newblock Concerning the elastico-plastic regime of filtration.
\newblock {\em Izv. Akad. Nauk SSSR, OTN}, (2):5--13, 1955.

\bibitem{barenblatt1969self}
GI~Barenblatt and GI~Sivashinskii.
\newblock Self-similar solutions of the second kind in nonlinear filtration:
  Pmm vol. 33, n 5, 1969, pp. 861--870.
\newblock {\em Journal of Applied Mathematics and Mechanics}, 33(5):836--845,
  1969.

\bibitem{barenblatt1972self}
GI~Barenblatt and Ya~B Zel'Dovich.
\newblock Self-similar solutions as intermediate asymptotics.
\newblock {\em Annual Review of Fluid Mechanics}, 4(1):285--312, 1972.

\bibitem{barenblatt1989theory}
Grigorij~I Barenblatt, Vladimir~Mordukhovich Entov, and
  Viktor~Mikha{\u\i}lovich Ryzhik.
\newblock Theory of fluid flows through natural rocks.
\newblock 1989.

\bibitem{barenblatt1996scaling}
Grigory~Isaakovich Barenblatt.
\newblock {\em Scaling, self-similarity, and intermediate asymptotics:
  dimensional analysis and intermediate asymptotics}, volume~14.
\newblock Cambridge University Press, 1996.

\bibitem{luis1995fully}
Luis~A. Caffarelli and Xavier Cabr{\'e}.
\newblock {\em Fully nonlinear elliptic equations}, volume~43.
\newblock American Mathematical Soc., 1995.

\bibitem{caffarelli2008counterexample}
Luis~A Caffarelli and Ulisse Stefanelli.
\newblock A counterexample to $\mathcal{C}^{2,1}$ regularity for parabolic
  fully nonlinear equations.
\newblock {\em Communications in Partial Differential Equations},
  33(7):1216--1234, 2008.

\bibitem{eggers2009role}
Jens Eggers and Marco~A Fontelos.
\newblock The role of self-similarity in singularities of partial differential
  equations.
\newblock {\em Nonlinearity}, 22(1):R1, 2009.

\bibitem{goldenfeld1990anomalous}
Nigel Goldenfeld, Olivier Martin, Y~Oono, and Fong Liu.
\newblock Anomalous dimensions and the renormalization group in a nonlinear
  diffusion process.
\newblock {\em Physical review letters}, 64(12):1361, 1990.

\bibitem{hu2012explicit}
Mingshang Hu.
\newblock Explicit solutions of the {G}-heat equation for a class of initial
  conditions.
\newblock {\em Nonlinear Analysis: Theory, Methods \& Applications},
  75(18):6588--6595, 2012.

\bibitem{HUANG2012LARGE}
Yong Huang and Juan~L V{\'a}zquez.
\newblock Large-time geometrical properties of solutions of the {B}arenblatt
  equation of elasto-plastic filtration.
\newblock {\em Journal of Differential Equations}, 252(7):4229--4242, 2012.

\bibitem{kamin1991barenblatt}
Shoshana Kamin, Lambertus~A Peletier, and Juan~Luis Vazquez.
\newblock On the {B}arenblatt equation of elastoplastic filtration.
\newblock {\em Indiana Univ. Math. J}, 40(4):1333--1362, 1991.

\bibitem{Kerchman1971self}
VI~Kerchman.
\newblock On self-similar solutions of the second kind in the theory of
  unsteady filtration: Pmm vol. 35, n 1, 1971, pp. 189--192.
\newblock {\em Journal of Applied Mathematics and Mechanics}, 35(1):158--162,
  1971.

\bibitem{khuzhaerov2007inverse}
B~Kh Khuzhaerov and {\'E}~Ch Kholiyarov.
\newblock Inverse problems of elastoplastic filtration of liquid in a porous
  medium.
\newblock {\em Journal of Engineering Physics and Thermophysics},
  80(3):517--525, 2007.

\bibitem{kroger1994regularity}
Pawel Kr{\"o}ger.
\newblock Regularity conditions on parabolic measures.
\newblock {\em Arkiv f{\"o}r Matematik}, 32(2):373--391, 1994.

\bibitem{olver1974introduction}
FWJ Olver.
\newblock Introduction to asymptotics and special functions.
\newblock {\em Academic Press}, 1974.

\bibitem{peng2007g}
Shige Peng.
\newblock G-expectation, {G}-brownian motion and related stochastic calculus of
  {I}t{\^o} type.
\newblock In {\em Stochastic analysis and applications}, pages 541--567.
  Springer, 2007.

\bibitem{peng2008multi}
Shige Peng.
\newblock Multi-dimensional {G}-brownian motion and related stochastic calculus
  under {G}-expectation.
\newblock {\em Stochastic Processes and their Applications},
  118(12):2223--2253, 2008.

\end{thebibliography}

\end{document}